\documentclass{amsart}
\usepackage{amsmath}
\usepackage{amsfonts}
\usepackage{graphics}
\usepackage{epsfig}
\usepackage{amssymb}
\usepackage{amscd}

\newtheorem{thm}{Theorem}[section]

\theoremstyle{definition}

\newtheorem*{ack}{Acknowledgement}
\theoremstyle{remark}
\newtheorem{rem}{Remark}

\theoremstyle{definition}

\numberwithin{equation}{section}
\numberwithin{figure}{section}

\def\Hom{{\text{\rm{Hom}}}}
\def\End{{\text{\rm{End}}}}

\def\tr{{\text{\rm{tr}}}}
\def\rchi{{\hbox{\raise1.5pt\hbox{$\chi$}}}}
\def\Aut{{\text{\rm{Aut}}}}

\def\isom{\cong}
\def\tensor{\otimes}
\def\dsum{\oplus}

\def\reg{{\text{\rm{reg}}}}
\def\Ghat{\hat{G}}
\def\lam{\lambda}
\def\vol{{\text{\rm{vol}}}}

\def\rank{{\text{\rm{rank}}\,}}
\def\Sym{{\text{\rm{Sym}}}}
\def\Spec{{\text{\rm{Spec}}}}
\def\Jac{{\text{\rm{Jac}}}}
\def\Pic{{\text{\rm{Pic}}}}
\def\Prym{{\text{\rm{Prym}}}}
\def\Nm{{\text{\rm{Nm}}}}

\begin{document}
\large

\title[Geometry of Character Varieties]{Geometry of 
Character Varieties of Surface Groups}
\author[Motohico Mulase]{Motohico Mulase$^*$}  
\address{
Department of Mathematics\\
University of California\\
Davis, CA 95616--8633}
\email{mulase@math.ucdavis.edu}
\thanks{$^*$Research supported by NSF grant DMS-0406077 and UC Davis.}
\begin{abstract}
This article is based on a talk delivered at the 
RIMS--OCAMI Joint International Conference on
\emph{Geometry Related to Integrable Systems}
in September, 2007. Its aim is to review a recent
progress in the Hitchin integrable systems and character varieties 
of the fundamental groups of Riemann surfaces.
A survey on geometric aspects of these character varieties
is also provided as we develop the exposition from a simple
case to more elaborate cases.
\end{abstract}

\maketitle

\allowdisplaybreaks

\tableofcontents

\section{Introduction}
\label{sect:intro}

The character varieties we  consider in this article
are the set of equivalence classes
$$
\Hom(\pi_1(\Sigma_g), G)/G
$$
of representations of a surface group $\pi_1(\Sigma_g)$
into another group $G$.
Here $\Sigma_g$ is a closed oriented surface of
genus $g$, which is   assumed to be
$g\ge 2$ most of the time. The action of $G$ on the
space of homomorphisms is through the conjugation action.
Since this action has fixed points, the quotient requires a special
treatment to make it a reasonable space. 
Despite the simple appearance of the space, it has an essential
 connection
to many other subjects in mathematics (\cite{AB, BNR, FQ,
FS, Fukaya, HN, H1, H2, JK, LM1, MFK, NS, W1991, Z}), and the list is 
steadily growing (\cite{DP, Frenkel, HRVK, HT, KW, MY2}). 
Our subject thus
 provides an ideal window 
to observe the scenery of a good part of recent developments in mathematics
and mathematical physics.

Each section of this article is devoted to a specific
type of character
varieties and a particular group $G$.
 We start with a finite group in Section~\ref{sect:finite}.
Already in this case one can appreciate the interplay between
the character variety and the theory of irreducible representations
of a finite group. In Sections \ref{sect:bundle}
and \ref{sect:twisted} we consider the  case $G=U_n$. 
We review the discovery of the relation to two-dimensional
Yang-Mills theory and symplectic geometry
due to Atiyah and Bott \cite{AB}. It forms the
turning point of the modern developments on character varieties.
We then turn our attention to the case $G=GL_n(\mathbb{C})$
in Sections \ref{sect:gln}  and \ref{sect:Hitchin}.
Here the key ideas
we review are due to Hitchin \cite{H1, H2}. 
In these seminal papers
Hitchin has suggested the subject's 
possible relations to four-dimensional 
Yang-Mills theory and the Langlands duality. These connections
are materialized recently by Hausel and Thaddeus \cite{HT},
Donagi and Pantev \cite{DP}, 
Kapustin and Witten \cite{KW},
 and many others. Section~\ref{sect:quotient}
motivates some of  these developments from our study \cite{HM} on 
the Hitchin integrable systems.

\section{Character varieties of finite groups and representation theory}
\label{sect:finite}

The simplest example of character varieties occurs when
$G$ is a finite group. The ``variety''  is a finite
set, and the only interesting
invariant is its cardinality. Here the \emph{reasonable} 
quotient  $\Hom(\pi_1(\Sigma_g), G)/G$
is \emph{not} the orbit space. A good theory
exists only for the \emph{virtual} quotient, which takes into
account the information of isotropy subgroups, exactly as
we do when we consider orbifolds.
\begin{thm}[Counting formula]
The classical counting formula gives
\begin{equation}
\label{eq:finite}
\frac{|\Hom(\pi_1(\Sigma_g),G)|}{|G|}
=
\sum_{\lam\in\Ghat} \left(\frac{\dim\lam}{|G|}\right)^
{\rchi(\Sigma_g)} ,
\end{equation}
where $\Ghat$ is the set of irreducible representations
of $G$, $\dim\lam$ is the dimension of the irreducible representation
$\lam\in\Ghat$, and $\rchi(\Sigma_g) = 2-2g $ is the Euler 
characteristic of the surface. 
\end{thm}
When $g=0$, the above formula
reduces to a well-known formula in representation theory:
\begin{equation}
\label{eq:order}
|G| = \sum_{\lam\in\Ghat} (\dim\lam)^2 .
\end{equation}

\begin{rem}
The formula for $g=1$ is known to Frobenius
\cite{F}. Burnside asks a related question as 
an exercise of his textbook \cite{Burnside}.
 In the late 20th century,
the formula was rediscovered by Witten \cite{W1991}
using quantum Yang-Mills theory in two dimensions,
and by Freed and Quinn \cite{FQ} using quantum
Chern-Simons gauge theory with the finite group $G$ as
its gauge group. \end{rem}

\begin{rem}
Since 't\,Hooft
\cite{'tHooft}  we know that a matrix integral
admits a ribbon graph expansion, using the Feynman diagram
technique \cite{Feynman}. In \cite{MY2} we 
ask what types of integrals admit a ribbon graph expansion.
Our answer is that an integral over a von Neumann algebra
admits such an expansion. We find in \cite{MY1, MY2} that
 when we apply
a formula of \cite{MY2} to the complex group 
algebra $\mathbb{C}[G]$, the counting formula
(\ref{eq:finite}) for 
all values of  $g$ automatically follows. The key fact is the algebra
decomposition
\begin{equation}
\label{eq:algebraiso}
\mathbb{C}[G] \isom \bigoplus_{\lam\in\Ghat} \End(\lam).
\end{equation}
The integral over the group algebra
then decomposes
into the product of matrix integrals over each simple factor $\End(\lam)$,
which we know how to calculate by 't Hooft's method.
Although  (\ref{eq:finite}) looks like a generalization of (\ref{eq:order}),
these formulas actually contain the same amount of information
because they are direct consequences of the decomposition 
(\ref{eq:algebraiso}). 
\end{rem}

\begin{rem}
We also note that there are corresponding 
formulas for  closed \emph{non-orientable}
surfaces \cite{MY1, MY2}. Intriguingly, the formula for non-orientable
surfaces are studied in its full generality, 
though without any mention on its geometric significance, in
a classical paper by Frobenius and Schur  \cite{FS}. 
The Frobenius-Schur theory automatically appears in 
the generalized matrix integral over the \emph{real}
group algebra $\mathbb{R}[G]$ (see \cite{MY1}). 
\end{rem}

Of course (\ref{eq:finite}) has an elementary proof,
without appealing to quantum field theories or matrix integrals.
We record  it
here only assuming a minimal 
background of representation theory that can be
found, for example, in Serre's textbook \cite{Serre}. 

The fundamental group of a compact oriented surface of genus
$g$ is generated by $2g$ generators with one relator:
$$
\pi_1(\Sigma_g) = \langle a_1,b_1,\dots,a_g,b_g\;|\;
[a_1,b_1]\cdots [a_g,b_g]=1\rangle ,
$$
where $[a,b] = aba^{-1}b^{-1}$. 
Since 
\begin{equation}
\label{eq:homspace}
\Hom(\pi_1(\Sigma_g),G)
= \{(s_1,t_1,\dots,s_g,t_g)\in G^{2g}\; \big| \;
[s_1,t_1]\cdots [s_g,t_g] = 1\} ,
\end{equation}
the counting problem reduces to evaluating an \emph{integral}
\begin{equation}
\label{eq:integral}
\big|\Hom(\pi_1(\Sigma_g),G)\big|
= \int_{G^{2g}}
\delta([s_1,t_1]\cdots [s_g,t_g] )ds_1dt_1\cdots ds_gdt_g .
\end{equation}
Here the left hand side is the volume of the character variety
that is defined by an invariant measure $ds$ on the group $G$.
For the case of a finite group, 
the volume is simply the cardinality, and the integral is
the sum over $G^{2g}$.
The $\delta$-function on $G$ is given by the normalized 
character of the regular representation
\begin{equation}
\label{eq:delta}
\delta(x) = \frac{1}{|G|}\rchi_\reg (x) =  \sum_{\lam\in\Ghat}
\frac{\dim \lam}{|G|}\cdot \rchi_\lam(x) .
\end{equation}
To compute the integral (\ref{eq:integral}),
let us first identify the complex group algebra 
$$
\mathbb{C}[G] = \bigg\{x=\sum_{\gamma\in G} x(\gamma)\cdot
\gamma\, \bigg|  \,x(\gamma)\in\mathbb{C}\bigg\}
$$
of a finite group $G$ 
with the vector space $F(G)$ of functions on $G$. 
In this way we can reduce the complexity of 
the commutator produce in (\ref{eq:homspace})
into simpler pieces. The
\emph{convolution product} of  two functions
$x(\gamma)$ and $y(\gamma)$ is defined by
$$
(x*y)(w) \overset{\text{def}}{=} \sum_{\gamma\in G}x(w\gamma^{-1})y(\gamma)\;,
$$
which makes $(F(G),*)$ an algebra 
 isomorphic to the group algebra. 
 In this identification,
the set of class functions $C\!F(G)$ corresponds to the center 
$Z\mathbb{C}[G]$ of 
$\mathbb{C}[G]$. 
According to the decomposition of this algebra into simple factors 
(\ref{eq:algebraiso}),
 we have an algebra 
isomorphism
$$
Z\mathbb{C}[G] = \bigoplus_{\lambda\in\hat{G}}
\mathbb{C}\;,
$$
where each factor $\mathbb{C}$ is the center of
$\End{\lambda}$. The projection to each factor is
given by 
$$
pr_\lambda : Z\mathbb{C}[G]\owns
x = \sum_{\gamma\in G} x(\gamma) \cdot \gamma
\longmapsto 
pr_\lambda(x)\overset{\text{def}}{=}\frac{1}{\dim \lambda}
\sum_{\gamma\in G} x(\gamma) \rchi_\lambda(\gamma)
\in\mathbb{C}\;,
$$
where $\rchi_\lam$ is the character of $\lam\in\Ghat$.
Following Serre \cite{Serre}, let
\begin{equation}
\label{eq:proj}
p_\lambda \overset{\text{def}}{=} \frac{\dim \lambda}{|G|}
\sum_{\gamma\in G} \rchi_\lambda (\gamma^{-1})\cdot
\gamma \in Z\mathbb{C}[G], \qquad \lam\in\Ghat ,
\end{equation}
 be a linear bases for $Z\mathbb{C}[G]$.
It follows from Schur's orthogonality of the irreducible characters that
$pr_\lambda(p_\mu) = \delta_{\lambda\mu}$. 
Consequently, we have
$p_\lambda p_\mu = \delta_{\lambda\mu} p_\lambda$, or 
equivalently, 
\begin{equation*}
\begin{split}
&\frac{\dim \lambda}{|G|}
\sum_{s\in G} \rchi_\lambda (s^{-1})\cdot
s \cdot 
\frac{\dim \mu}{|G|}
\sum_{t\in G} \rchi_\mu (t^{-1})\cdot t\\
&=\frac{\dim \lambda\cdot \dim \mu}{|G|^2}
\sum_{w\in G}\left(
\sum_{t\in G} \rchi_\lambda ((wt^{-1})^{-1})
\rchi_\mu (t^{-1})\right)
\cdot w\\
&=\delta_{\lambda\mu}
\frac{\dim \lambda}{|G|}
\sum_{w\in G} \rchi_\lambda (w^{-1})\cdot w\;.
\end{split}
\end{equation*}
We thus obtain
\begin{equation}
\label{eq:charconv}
\rchi_\lambda *\rchi_\mu = 
\frac{|G|}{\dim \mu}\delta_{\lambda\mu}
\rchi_\lambda\;.
\end{equation}
We now turn to the counting formula. 
Let 
\begin{equation}
\label{eq:fg}
f_g (w) \overset{\text{def}}{=}\big|
\{(s_1,t_1,s_2,t_2,\dots,s_g,t_g)
\in G^{2g}\;|
\; [s_1,t_1]
\cdots [s_g,t_g] =
w\}\big| .
\end{equation}
This is a class function and satisfies
$f_g(w) = f_g(w^{-1})$. 
{F}rom the definition, it is obvious that
$
f_{g_1 + g_2} = f_{g_1} * f_{g_2}
$.
Therefore, 
\begin{equation}
\label{eq:fgfromf1}
f_g = \overset{g\text{-times}}{\overbrace{f_1 *\cdots * f_1}}
\;.
\end{equation}
Finding $f_1$
is Exercise~7.68 of Stanley's textbook \cite{Stanley},
and the answer is in Frobenius \cite{F}.
From Schur's lemma,
\begin{equation}
\label{eq:irrepsum}
 \sum_{s\in G} \rho_\lam(s\cdot t\cdot s^{-1})
\end{equation}
is central as an element of $\End(\lam)$, where $\rho_\lam$ is the irreducible representation corresponding
to $\lam\in\Ghat$. This is because 
(\ref{eq:irrepsum}) commutes
with $\rho_\lam(w)$ for every $w\in G$. Hence we have
$$
\sum_{s\in G} \rho_\lam(s\cdot t\cdot s^{-1})
= \sum_{s\in G} \frac{\rchi_\lam(s\cdot t\cdot s^{-1})}{\dim\lam}
=\frac{|G|}{\dim\lam}\rchi_\lam(t) ,
$$
noticing that the character $\rchi_\lam$ is the trace of $\rho_\lam$. 
Therefore, 
\begin{equation*}
\begin{aligned}
\dim\lam \sum_{s\in G} \rho_\lam (s\cdot t\cdot s^{-1}\cdot t^{-1}w^{-1})&=
\dim\lam \sum_{s\in G} \rho_\lam (s\cdot t\cdot s^{-1})\cdot \rho_\lam( t^{-1}w^{-1})\\
&=|G|\cdot  \rchi_\lam(t)\cdot \rho_\lam(t^{-1}w^{-1})\;.
\end{aligned}
\end{equation*}
Taking trace and summing in $t\in G$ of the above equality, we obtain
$$
\frac{\dim\lam}{|G|}\sum_{s,t\in G}  \rchi_\lam(sts^{-1}t^{-1}w^{-1})
=\sum_{t\in G}\rchi_\lam(t)\rchi_\lam(t^{-1}w^{-1})=
( \rchi_\lam *\rchi_\lam)(w^{-1}) = \frac{|G|}{\dim\lam}\cdot \rchi_\lam(w^{-1}).
$$
Switching to the $\delta$-function of
(\ref{eq:delta}), we find
\begin{equation}
\label{eq:f1}
f_1(w)=\int_{G^2} \delta([s,t]w^{-1})dsdt = \sum_{\lam\in\Ghat}
\frac{|G|}{\dim\lam}\cdot \rchi_\lam(w^{-1})
= \sum_{\lam\in\Ghat}
\frac{|G|}{\dim\lam}\cdot \rchi_\lam(w) .
\end{equation}
Note that we can interchange $w$ and $w^{-1}$, since $f_g$ is
integer valued and is invariant under complex conjugation.
From (\ref{eq:charconv}),  (\ref{eq:fgfromf1}) and (\ref{eq:f1}), we obtain
\begin{thm}[Counting formula for twisted case]
For every $g\ge 1$ and $w\in G$ let 
$$
f_g (w) {=}\big|
\{(s_1,t_1,s_2,t_2,\dots,s_g,t_g)
\in G^{2g}\;|
\; [s_1,t_1]
\cdots [s_g,t_g] =
w\}\big| .
$$
Then we have a character expansion formula
\begin{equation}
\label{eq:charexpansion}
f_g(w) = f_g(w^{-1}) = \sum_{\lambda\in \hat{G}} \left(
\frac{|G|}{\dim\lambda}\right)^{2g-1}
\cdot \rchi_\lambda (w)\;.
\end{equation}
\end{thm}
The counting formula (\ref{eq:finite}) is a special
case for $w=1$.

\section{Character varieties of $U_n$ as moduli spaces of stable vector bundles}
\label{sect:bundle}
 
The next natural case of character varieties is for  a compact
Lie group $G$,  in particular, $G=U_n$.
The issue of taking the quotient $\Hom(\pi_1(\Sigma_g),U_n)/U_n$
is much more serious than the finite group case, due to the fact that
the trivial representation of $\pi_1(\Sigma_g)$ into $U_n$ is a 
fixed point of the conjugation action. Consequently, the quotient
space does not have a good manifold structure
 at the trivial representation. One way 
to avoid this and other quotient difficulties is to 
 restrict our consideration to  
\emph{irreducible} unitary representations
\begin{equation}
\label{eq:Unirred}
\Hom^{\text{irred}}(\pi_1(\Sigma_g),U_n)/U_n .
\end{equation}
From now on we assume $g\ge 2$.
This time the quotient is well-defined as a real analytic space
with some minor singularities. 
According to Narasimhan and Seshadri \cite{NS}, (\ref{eq:Unirred}) is
diffeomorphic to the moduli space, denoted here by
$\mathcal{U}_C(n,0)$, 
of stable holomorphic 
vector bundles of rank $n$ and degree $0$ on a smooth algebraic curve $C$ of genus $g$.
A holomorphic vector bundle $E$ on $C$ is said to be \emph{semistable} if
\begin{equation}
\label{eq:semistable}
\frac{\deg F}{\rank F}\le \frac{\deg E}{\rank E}
\end{equation}
for every holomorphic proper vector subbundle $F\subset E$, and
\emph{stable} if the strict inequality holds.
If the rank and the degree 
are relatively prime, then the equality cannot
hold  in (\ref{eq:semistable}),
hence every semistable vector bundle is automatically stable.
The topological structure of a vector bundle $E$ on $\Sigma_g$ is 
determined by its rank and the degree. 
From the expression (\ref{eq:Unirred}) it is clear that the differentiable
structure of $\mathcal{U}_\mathcal{C}(n,0)$ does not depend on which 
complex structure we give on $\Sigma_g$. 

As explained in the newest addition to Mumford's textbook
\cite{MFK} by Kirwan, moduli theory of
stable objects can also be understood in terms of 
the \emph{symplectic quotient} of the space of differentiable
connections on $C$ with values in $U_n$ by the group of gauge
transformations. Let $E$ be a topologically trivial
differentiable $U_n$-vector bundle on
$\Sigma_g$, and $\mathcal{A}(\Sigma_g,U_n)$ the space of 
differentiable connections in $E$. We denote by $ad(E)$ the 
associated adjoint $u_n$-bundle on $\Sigma_g$.
Since the tangent space to the space of $U_n$-connections is 
the space of sections  $\Gamma(\Sigma_g, ad(E)\tensor \Lambda^1(\Sigma_g))$, we can define a gauge invariant symplectic form
\begin{equation}
\label{eq:symplectic}
\omega(\alpha, \beta) = \frac{1}{8\pi^2}\int_C \tr (\alpha\wedge\beta),
\qquad \alpha,\beta\in 
\Gamma(\Sigma_g, ad(E)\tensor \Lambda^1(\Sigma_g))
\end{equation}
on the space of 
$U_n$-connections on $\Sigma_g$. The Lie algebra of the group 
$\mathcal{G}(\Sigma_g,U_n)$ of 
gauge transformations is the space of global sections of $ad(E)$, 
hence its dual is $\Gamma(\Sigma_g, ad(E)\tensor \Lambda^2(\Sigma_g))$. The moment map of the 
$\mathcal{G}(\Sigma_g,U_n)$-action
on the space of connections is then given by the curvature map
\begin{equation}
\label{eq:moment}
\mu_\Sigma: \mathcal{A}(\Sigma_g,U_n)\owns A\longmapsto
F_A = dA + A\wedge A\in 
\Gamma(\Sigma_g, ad(E)\tensor \Lambda  ^2(\Sigma_g)).
\end{equation}
If we choose $0\in \Gamma(\Sigma_g, ad(E)\tensor \Lambda^2(\Sigma_g))$ as the reference value of the moment map, then 
the symplectic quotient 
$$
\mathcal{A}(\Sigma_g,U_n)/\!\!/\mathcal{G}(\Sigma_g,U_n)
= \mu_\Sigma ^{-1}(0)/\mathcal{G}(\Sigma_g,U_n)
=\Hom(\pi_1(\Sigma_g),U_n)/U_n
$$
gives the moduli space of flat $U_n$-connections on 
$\Sigma_g$. This correspondence is also known as 
the Riemann-Hilbert correspondence.

If the structure of a compact Riemann surface $C$ is 
chosen on $\Sigma_g$, 
then a connection in a differentiable vector bundle $E$ on $C$ 
defines a holomorphic structure in $E$. This process
goes as follows. First we note that  there are no type $(0,2)$-forms 
on $C$. Therefore, the $(0,1)$-part of the connection is
always integrable. We can then define a differentiable
section of $E$ to be \emph{holomorphic}
if it is annihilated by 
the $(0,1)$-part of the covanriant derivative.
If the connection $A$ is unitary, then it is 
uniquely determined by it's $(0,1)$-part. The information of
$A$ is thus encoded in the complex structure it defines on $E$. 
In particular, the moduli space of flat unitary connections
modulo gauge equivalence becomes
the moduli space of holomorphic vector bundles of degree $0$. 
The stability condition of a holomorphic vector bundle
is equivalent to requiring that the corresponding flat connection
is irreducible. This in turn corresponds to
irreducibility of the unitary representation of
$\pi_1(C)$. Since the  curvature $F_A$
receives a topological constraint, the moment map
(\ref{eq:moment})
  cannot take an arbitrary value of
$\Gamma(\Sigma_g, ad(E)\tensor \Lambda^2(\Sigma_g))$. In particular, 
$0$
is a critical value of the moment map $\mu_\Sigma$, and hence
the symplectic quotient is
 singular. 
 
Although we have this issue of singularities, the above discussion 
 shows that the $U_n$-character variety outside its
singularities has a
natural symplectic structure coming from (\ref{eq:symplectic})
and the process of symplectic quotient, 
and a complex structure as the moduli space of holomorphic vector
bundles if a complex structure is chosen on $\Sigma_g$. 
The symplectic and complex structures are compatible,
so outside  the singularities the character variety is a
complex K\"ahler manifold. 
Consequently, its dimension 
should be even. Actually, we can 
compute the dimension directly from
(\ref{eq:homspace}). Noticing that $\det [s,t]= 1$ and that the
center of $U_n$ acts trivially via conjugation, we have 
\begin{equation}
\label{eq:dimension}
\dim_\mathbb{R} \Hom(\pi_1(\Sigma_g),U_n)/U_n
= n^2 (2g-2) + 2 = 2 (n^2 (g-1) +1).
\end{equation}

All the considerations become much simpler when the 
group is $G=U_1$. The condition of (\ref{eq:homspace}) is
vacuous and the  character variety is simply a $2g$-dimensional
real torus
$$
\Hom(\pi_1(\Sigma_g), U_1) = \Hom(H_1(\Sigma_g,\mathbb{Z}), U_1)
= (U_1)^{2g}.
$$
If a complex structure $C$ is chosen on $\Sigma_g$, then 
the complex line bundle arising from a representation
of $\pi_1(\Sigma_g)$ acquires a holomorphic structure,
and the character variety becomes the Jacobian:
$$
\Hom(\pi_1(C), U_1) \isom \Jac(C) = \Pic^0(C).
$$

\section{Twisted character varieties of $U_n$}
\label{sect:twisted}

To study moduli spaces of holomorphic vector bundles
on a Riemann surface that are not topologically trivial,
we need to consider a variant of character varieties.
Let $E$ now be a topological vector bundle of rank $n$
and degree $d\ne 0$ on $C=\Sigma_g$. This time it admits no flat connections,
because the degree of $E$ is determined by its connection through 
the Chern-Weil formula:
$$
\deg E = c_1(E) = -\frac{1}{2\pi i} \int_C \tr(F_A).
$$
The symplectic quotient of the space of connections in $E$ 
requires a point in the dual Lie algebra 
$
F_A\in \Gamma(\Sigma_g, ad(E)\tensor \Lambda^2(\Sigma_g))
$
that is fixed 
under the coadjoint action of $\mathcal{G}(\Sigma_g,U_n)$.
Obviously, $F_A$ is coadjoint invariant if it takes central values.
A unitary connection $A$ in $E$ is said to be
\emph{projectively flat} if its curvature $F_A$ is central. 
Narasimhan-Seshadri \cite{NS} again tells us that 
the moduli space $\mathcal{U}_C(n,d)$ of stable holomorphic
vector bundles on $C$ of rank $n$ and degree $d$ is
diffeomorphic to the space of gauge equivalent classes
of irreducible projectively flat connections. 

Among the projectively flat connections, there is a particularly
natural class. Since the curvature $F_A$ of a connection $A$ is
a 2-form, we cannot talk about $F_A$ being a constant. But
if we apply the Hodge $*$-operator, then the covariant constant 
condition
\begin{equation}
\label{eq:YM}
d_A * F_A = 0
\end{equation}
makes sense. This is exactly the two-dimensional
\emph{Yang-Mills equation}
studied by Atiyah and Bott in \cite{AB}. 
A projectly flat solution $A$ of the Yang-Mills equation has
its curvature given by 
\begin{equation}
\label{eq:curvature}
F_A = - \; \frac{2\pi i d}{n} I_n \cdot \vol_C ,
\end{equation} 
where $\vol_C$ is the normalized volume form of $C$ 
with total volume $1$. The \emph{holonomy
group} of a connection at a point $p\in C$ is generated by 
parallel transports along every closed loop that starts at $p$. 
The Lie algebra of the holonomy group is the Lie subalgebra
of $u_n$ in which the curvature form $F_A$ takes values.
For a projectively flat connection, the holonomy group is
the center $U_1$ of $U_n$. Certainly, the Lie algebra generated
by the value (\ref{eq:curvature}) is $\mathbb{R}$, and the
corresponding Lie group is $U_1$.

The Riemann-Hilbert correspondence
gives an identification
between a flat connection and a representation of $\pi_1(\Sigma_g)$
into $U_n$. What is a counterpart of the Riemann-Hilbert
correspondence for the case of a projectively flat connection?

When the curvature is non-zero, a parallel transport of 
a connection does not induce a representation
$\pi_1(\Sigma_g)\rightarrow U_n$ because it depends on
the choice of a loop.
The answer to the above question presented in \cite{AB}
is that \emph{a projective Yang-Mills connection corresponds to
a representation of a central extension of $\pi_1(\Sigma_g)$ into 
$U_n$.} In the following we examine this correspondence
for irreducible connections. 

We note that $\pi_1(\Sigma_g)$ has a universal central
extension
\begin{equation}
\label{eq:central}
\begin{CD}
1@>>>\mathbb{Z}@>>>\hat{\pi}_1(\Sigma_g)@>>>{\pi}_1(\Sigma_g)
@>>>1,
\end{CD}
\end{equation}
where the extended group is defined by
$$
\hat{\pi}_1(\Sigma_g) =\langle a_1,b_1,\dots,a_g,b_g,c\;|\;
[c,a_i]=[c,b_i]=1, [a_1,b_1]\cdots [a_g,b_g]=c\rangle ,
$$
and $\mathbb{Z}\owns k\longmapsto c^k\in \hat{\pi}_1(\Sigma_g)$
determines its center. The central extension we need is a \emph{Lie group}
$\hat{\pi}_1(\Sigma_g)_\mathbb{R}$
that contains a copy of $\mathbb{R}$ through
$
\mathbb{R}\owns r\longmapsto c^r\in \hat{\pi}_1(\Sigma_g)_\mathbb{R}
$,
and satisfies that
\begin{equation}
\label{eq:exact}
\begin{CD}
1@>>>\mathbb{R}@>>>\hat{\pi}_1(\Sigma_g)_\mathbb{R}@>>>{\pi}_1(\Sigma_g)
@>>>1.
\end{CD}
\end{equation}

\begin{thm}[Atiyah-Bott \cite{AB}] The
\textbf{twisted} character variety
\begin{equation}
\label{eq:twisted}
\Hom^{\emph{irred}}(\hat{\pi}_1(\Sigma_g)_\mathbb{R}, U_n)/U_n
\end{equation}
of irreducible representations
is identified with the space of irreducible unitary Yang-Mills 
connections in $E$ modulo gauge transformations. 
\end{thm}

Note that 
$\Hom(\hat{\pi}_1(\Sigma_g)_\mathbb{R}, U_n) =
\{ (s_1,t_1,\dots,s_g,t_g,\gamma)\in (U_n)^{2g+1}\;|\;
[\gamma,s_i] = [\gamma,t_i] = 1, \;[s_1,t_1]\cdots [s_g,t_g] = \gamma\}
$.
Since the commutator product is equated to $\gamma\in U_n$
which is not necessarily the identity, the name ``twisted''
is used in the literature.
If a representation $\hat{\pi}_1(\Sigma_g)_\mathbb{R}\rightarrow
U_n$ is irreducible, then $\gamma$ is a central element of $U_n$. 
Since $\det [s,t] = 1$, we conclude that 
\begin{equation}
\label{eq:gammavalue}
\gamma= \exp\left(\frac{2\pi i d}{n}\right)\cdot I_n
\end{equation}
for some integer $d$. Therefore, $\Hom^{\text{irred}}
(\hat{\pi}_1(\Sigma_g)_\mathbb{R}, U_n)$ consists
of $n$ disjoint pieces corresponding to the $n$
possible values for (\ref{eq:gammavalue}).

The construction of a Yang-Mills connection from
an irreducible representation 
$$
\rho\in \Hom^{\text{irred}}(\hat{\pi}_1(\Sigma_g)_\mathbb{R}, U_n)
$$ 
goes as follows. First we choose a connection $a$ in a 
complex line bundle $L$ on $\Sigma_g $ of degree $1$. 
The Yang-Mills equation for $a$ is simply the linear 
harmonic equation $d*da=0$ because
$U_1$ is Abelian. So let us choose a harmonic connection
$a$ with curvature 
\begin{equation}
\label{eq:Fa}
F_a = -{2\pi i}\cdot \vol_\Sigma .
\end{equation}
Let $h: \hat{\Sigma}_g\rightarrow \Sigma_g$ be the universal
covering of $\Sigma_g$. Then the pull-back line bundle $h^*L$
on $\hat{\Sigma}_g$, viewed as a fiber  bundle on $\Sigma_g$, 
has the structure group $U_1\times \pi_1(\Sigma_g)$.
Note that the exact sequence (\ref{eq:exact}) induces a
surjective homomorphism 
$$
f: \hat{\pi}_1(\Sigma_g)_\mathbb{R} \longrightarrow U_1\times 
\pi_1(\Sigma_g)
$$
by sending the central generator $c$ to a non-identity
element of $U_1$. We can thus construct a principal 
$\hat{\pi}_1(\Sigma_g)_\mathbb{R}$-bundle $P$ on $\Sigma_g$
from $L$, $h$, and $f$, in which the lift of $a$ now lives as
a Yang-Mills connection with the constant curvature
(\ref{eq:Fa}). Consider the principal $U_n$-bundle on
$\Sigma_g$ defined by $P\times_\rho U_n$, and its associated
rank $n$ vector bundle $E$ through the standard $n$-dimensional
representation of $U_n$ on $\mathbb{C}^n$. Let $A$ be the natural 
connection in $E$ arising from $a$. Then by functoriality
of the Yang-Mills equation, $A$ is automatically a
Yang-Mills connection in $E$. The holonomy of $A$ is the group 
generated by $\gamma = \rho(c)$ in $U_n$, which is central since
$\rho$ is irreducible. The value of the curvature $F_A$ of $A$ 
is quantized according to the topological type of $E$,
which is also determined by $\rho(c)\in U_n$. 

To show that every irreducible unitary Yang-Mills connection 
gives rise to a representation 
$$
\rho: \hat{\pi}_1(\Sigma_g)_\mathbb{R} \rightarrow U_n ,
$$
first we note that the same statement is true for $G=U_1$ and 
$G=SU_n$. Then we reduce the problem of construction
to the
hybrid of these two cases.
For $SU_n$, the vector bundle involved is trivial,
and an irreducible Yang-Mills connection is necessarily flat. 
Thus it gives rise to a representation of $\pi_1(\Sigma_g)$. 
For $U_1$, the group is Abelian and the question reduces to 
the standard homology theory. By pulling back
a unitary connection through the covering homomorphism
$$
U_1\times SU_n \longrightarrow U_n ,
$$
we can reduce the general case to the two special cases 
\cite{AB}.

An important fact is that if $\gamma$ of (\ref{eq:gammavalue})
is a primitive root of unity, i.e., $G.C.D.(n,d) = 1$, then
$\mathcal{U}_C(n,d)$ is  
 a non-singular projective algebraic variety. 
The smoothness is a consequence 
of the fact that such a $\gamma$ is a regular value of
the commutator product map
\begin{equation}
\label{eq:commutator}
\mu:  (U_n)^{2g}\owns (s_1,t_1,\dots,s_g,t_g)\longmapsto
[s_1,t_1]\cdots[s_g,t_g]\in SU_n ,
\end{equation}
and that the isotropy subgroup of the conjugation action of $U_n$
on $\mu^{-1}(\gamma)$ is always the central $U_1$. 
These statements are easily verified through direct calculations
(see for example \cite{HRVK}). 
Let us choose a point 
 $p = (s_1,t_1,\dots,s_g,t_g)\in \mu^{-1}(\gamma)$
in the inverse image of a primitive root of unity $\gamma$. 
The differential $d\mu_p$ of $\mu$ at $p$ is a linear map
between Lie algebras
$$
d\mu_p : (u_n)^{\dsum 2g} \longrightarrow su_n .
$$
Note that for $s\in U_n$ and $x\in u_n$, we have $ds(x) = x$.
Let us first consider the case $g=1$. We wish
to show that 
\begin{equation*}
\begin{aligned}
&d\mu_p (x,y)\\
 = &ds(x)\cdot ts^{-1}t^{-1} + s\cdot dt(y)\cdot s^{-1}
t^{-1} - sts^{-1}\cdot ds(x)\cdot s^{-1}t^{-1} 
-sts^{-1}t^{-1}\cdot dt(y)\cdot t^{-1}\\
= &xts^{-1}t^{-1} + sys^{-1}t^{-1} - sts^{-1}xs^{-1}t^{-1} - sts^{-1}t^{-1} yt^{-1} \\
= &\gamma (xs^{-1} - t x s^{-1} t^{-1}) + \gamma (syt^{-1}s^{-1}
- yt^{-1})
\end{aligned}
\end{equation*}
spans the entire Lie algebra $su_n$ as $(x,y)\in (u_n)^2$ varies.
In the above computation
 products and additions are calculated as $n\times n$
complex matrices, and we have used the commutation relation
$sts^{-1}t^{-1}  =\gamma$. Recall that $\tr(vw)$ defines a non-degenerate
bilinear form on $su_n$. Suppose now that
$\tr(w \cdot d\mu_p(x,y))=0$ for all $x,y\in u_n$. For $y=0$ it
follows that 
\begin{equation*}
\begin{aligned}
&\tr(xs^{-1}w) = \tr(txs^{-1}t^{-1}w) \qquad \text{for all } x\in u_n\\
\Longleftrightarrow \quad &
s^{-1}w = s^{-1}t^{-1}wt\\
\Longleftrightarrow \quad &w = t^{-1}wt.
\end{aligned}
\end{equation*}
Similarly, for $x=0$, we obtain $w=s^{-1}ws$. Therefore,
$w$ commutes with $s$ and $t$. We can then restrict
the relation $[s,t]=\gamma$ to any eigenspace of $w$ 
of dimension $m\le n$. The determinant condition $\det [s,t]=1$
yields $\gamma^m=1$. Hence $m=n$ because $\gamma$ is
primitive,  establishing
that $w$ is a scalar diagonal matrix. Since $w\in su_n$, we conclude
that $w=0$. 

For $g\ge 2$, we use the relation $[s_1,t_1]\cdots
[s_g,t_g]=\gamma$ to establish that any $w\in su_n$ that
satisfies $\tr(w\cdot d\mu_p(x_1,y_1,\dots,x_g,y_g))=0$
commutes with $s_1$ and $t_1$ when restricted to 
$x_i=y_i=0$ for $i>1$. We can then recursively show that
$w$ actually commutes with all $s_i$ and $t_i$. Restricting the
commutator product relation to any eigenspace of $w$
as above and using the fact that $\gamma$ is primitive,
we conclude that $w$ is central, and hence equal to $0\in su_n$. 
It follows that 
 $\gamma\in SU_n$ is 
a regular value of (\ref{eq:commutator}), and consequently
$\mu^{-1}(\gamma)$ is a non-singular manifold. 

Note that in the above argument we have also shown 
that the isotropy subgroup of $U_n$
acting on $\mu^{-1}(\gamma)$ through conjugation 
is the central $U_1$ at any point of $\mu^{-1}(\gamma)$.
Therefore, the quotient 
$$
\mu^{-1}(\gamma)/U_n = \mathcal{U}_C(n,d)
$$
is non-singular if $G.C.D.(n,d)=1$. 

The task of calculating the Poinar\'e polynomial 
of this non-singular compact complex algebraic manifold
is carried out by Harder-Narasimhan \cite{HN}, Atiyah-Bott
\cite{AB} and Zagier \cite{Z}. Harder and Narasimhan 
use Deligne's solution to the Weil conjecture (see for example
\cite{O}) as their tool
and study the moduli theory over the finite field $\mathbb{F}_q$
for all possible values of $q=p^e$. Atiyah and Bott use
2-dimensional Yang-Mills theory and equivariant Morse-Bott
theory to derive the topological structure of 
$\mathcal{U}_C(n,d)$. Both \cite{HN} and
\cite{AB} lead to a recursion formula for the 
Poincar\'e polynomials. Zagier \cite{Z} obtains a closed formula,
solving the recursion relation.

\section{Twisted character varieties of $GL_n(\mathbb{C})$}
\label{sect:gln}

Twisted character varieties
\begin{equation}
\label{eq:reductive}
\Hom(\hat{\pi}_1(\Sigma_g), G)/\!\!/G
\end{equation}
for a complex reductive group $G$ have received much attention 
in recent years from many different points of view 
\cite{DP, HRVK, HT, KW}. In this section we consider
the case $G=GL_n(\mathbb{C})$. 
The quotient (\ref{eq:reductive}) is a 
\emph{geometric invariant theory
quotient} of \cite{MFK}, due to the fact that $G$ is not compact. 
The categorical quotient contains the geometric quotient
$$
\Hom^{\text{irred}}(\hat{\pi}_1(\Sigma_g), GL_n(\mathbb{C}))
/GL_n(\mathbb{C}) .
$$
The argument of Section~\ref{sect:twisted} applies here
to show that the central generator $c\in \hat{\pi}_1(\Sigma_g)$
is mapped to a central element $\gamma\in GL_n(\mathbb{C})$,
which takes the same value as in (\ref{eq:gammavalue}).
Thus the character variety 
 consists of $n$ disjoint pieces,
and a component corresponding to a primitive $n$-th roots of unity 
is a non-singular affine algebraic subvariety of 
complex dimension $2(n^2(g-1) + 1)$ contained in 
$\mathbb{C}^{2gn^2}$. From now on we refer to this
non-singular piece at  a primitive $n$-th root of unity  $\gamma$ by 
\begin{equation}
\label{eq:X}
\mathcal{X}(\mathbb{C}) = 
\{\rho\in \Hom^{\text{irred}}(\hat{\pi}_1(\Sigma_g), GL_n(\mathbb{C})) \;|\; \rho(c)=\gamma\}/GL_n(\mathbb{C}).
\end{equation}

A surprising result recently obtained by Hausel, Rodriguez-Villegas
and Katz in \cite{HRVK} is the calculation 
of the mixed Hodge polynomial of this
character variety. Their key idea is Deligne's Hodge theory.
It states that the mixed Hodge 
polynomial of a complex algebraic variety 
$X(\mathbb{C})$ can be determined if one knows the cardinality
of the mod $q=p^e$ reduction $X(\mathbb{F}_q)$ of $X$ for every
prime $p$ (or most of them at least)
 and its power $e$. For the case of the character variety for $GL_n(\mathbb{C})$,
since its defining equation 
$$
[s_1,t_1]\cdots [s_g,t_g] = \gamma
$$
 is a set of polynomial
equations  defined over $\mathbb{Z}[\gamma]$
among the entries of the matrices, 
the mod $q$ reduction is given by $\mathcal{X}(\mathbb{F}_q)$
if $p$ is not a factor of $n$. 
Now the group $GL_n(\mathbb{F}_q)$ is finite,
so the cardinality of the character
 variety  is readily available from (\ref{eq:charexpansion})!

Since $U_n$ is the compact real form of $GL_n(\mathbb{C})$, 
the compact complex manifold $\mathcal{U}_C(n,d)$ is 
contained as the real part of  $\mathcal{X}(\mathbb{C})$
if $\gamma = \exp(2\pi i d/n)$ and $G.C.D.(n,d) = 1$. 
What is the relation between the complex structure of 
$\mathcal{X}(\mathbb{C})$ naturally arising from 
$GL_n(\mathbb{C})$ and that of $\mathcal{U}_C(n,d)$
coming from a complex structure $C$ on the surface
$\Sigma_g$? This question is addressed
in Section~\ref{sect:quotient}. 

If we view the non-singluar compact complex projective
algebraic variety $\mathcal{U}_C(n,d)$ as a real analytic
Riemannian manifold whose metric is determined by the
K\"ahler structure, then its \emph{complexification} 
is the total space of the cotangent bundle 
$T^*\mathcal{U}_C(n,d)$. This is because the canonical
symplectic form on $T^*\mathcal{U}_C(n,d)$ and the
Riemannian metric induced from $\mathcal{U}_C(n,d)$ together
determine the unique almost complex structure on the
cotangent bundle which is integrable. Since $\mathcal{X}(\mathbb{C})$
is a complexification of $\mathcal{U}_C(n,d)$, it contains
this cotangent bundle as a complex submanifold:
\begin{equation}
\label{eq:cotangentinX}
T^*\mathcal{U}_C(n,d)\subset \mathcal{X}(\mathbb{C}).
\end{equation}
Of course this embedding is \emph{never} a holomorphic map
with respect 
to the complex structure of $\mathcal{U}_C(n,d)$. 
So far we have noticed that there are at least \emph{two} different
complex structures in $T^*\mathcal{U}_C(n,d)$. One is what
we have just described as a complex submanifold of 
$\mathcal{X}(\mathbb{C})$, which we denote by 
$J$, and the other comes from
the  cotangent bundle of the complex manifold
$\mathcal{U}_C(n,d)$ denoted by $I$. These complex 
structures are indeed different, since an affine manifold
$\mathcal{X}(\mathbb{C})$ cannot contain a compact
complex manifold $\mathcal{U}_C(n,d)$ in it.

In this section we study the structure of $\mathcal{X}(\mathbb{C})$
from the point of view of $2$-dimensional
Yang-Mills theory following 
Hitchin \cite{H1}. Let us consider a topological complex vector
bundle $E$ of rank $n$ and degree $d$
on a Riemann surface $C$ of genus $g$,
and a complex connection $A_\mathbb{C}$ in $E$ with values in 
$gl_n(\mathbb{C})$. We choose a Hermitian fiber metric
in $E$ and reduce the structure group to $U_n$. 
The skew-Hermitian part $A$
of $A_\mathbb{C}$ is a unitary connection which is 
well-defined under the unitary gauge transformation 
$\mathcal{G}(C,U_n)$, though the whole gauge transformation
$\mathcal{G}(C,GL_n(\mathbb{C}))$ does not preserve the
skew-Hermitian part. Note that the action of $\mathcal{G}(C,U_n)$
on the Hermitian part  of $A_\mathbb{C}$ is a linear
transformation because a unitary gauge transformation of the $0$ connection
is skew-Hermitian. Therefore the Hermitian part 
$\Phi$ of $A_\mathbb{C}$
 can be identified as
a differential $1$-form on $C$ with values in $ad_\mathbb{C}(E)$, the $gl_n(\mathbb{C})$-bundle
associated to $ad(E)$:
$$
\Phi\in \Gamma(C,ad_\mathbb{C}(E)\tensor \Lambda^{1}(\Sigma_g)).
$$
 Using the complex coordinate on $C$, let $\phi$ be the  
 type $(1,0)$-part of $\Phi$: 
$$
\phi=\Phi^{(1,0)}\in \Gamma(C,ad_\mathbb{C}(E)\tensor \Lambda^{(1,0)}(C)).
$$
Here again $\phi$ is well-defined under
the unitary gauge transformation, and it uniquely determines 
$\Phi$ because of the Hermitian condition. 
In this way we obtain a $\mathcal{G}(C,U_n)$-space isomorphism
\begin{equation}
\label{eq:Aphi}
\mathcal{A}(C,GL_n(\mathbb{C})) \isom
\mathcal{A}(C,U_n)\times \Gamma(C,ad_\mathbb{C}(E)\tensor \Lambda^{(1,0)}(C)),
\end{equation}
which identifies $A_\mathbb{C}$ with the pair $(A,\phi)$
thus obtained. 
We will come back to the point
of the action of $\mathcal{G}(C,GL_n(\mathbb{C}))$
on these spaces a little later. 

Hitchin 
shows that the moment map on $\mathcal{A}(C,U_n)\times \Gamma(C,ad_\mathbb{C}(E)\tensor \Lambda^{(1,0)}(C))$ for the 
gauge group $\mathcal{G}(C,U_n)$-action is given by
$$
\mu_H : \mathcal{A}(C,U_n)\times \Gamma(C,ad_\mathbb{C}(E)\tensor \Lambda^{(1,0)}(C))
\owns (A,\phi)\longmapsto 
F_A + [\phi,\phi^*]\in 
\Gamma(C,ad(E)\tensor \Lambda^{(1,1)}(C)),
$$
where $F_A$ is the curvature form of $A$ and 
$[\phi,\phi^*] = \phi\wedge\phi^*+\phi^*\wedge\phi$
is an $ad(E)$-valued (i.e., a locally
\emph{skew}-Hermitian) $(1,1)$-form
on $C$. 
Although $\mathcal{A}(C,U_n)/\!\!/\mathcal{G}(C,U_n)$
is finite-dimensional, the symplectic quotient 
$\mu_H ^{-1}(0)/\mathcal{G}(C,U_n)$ is still infinite-dimensional
due to the second factor
$\Gamma(C,ad_\mathbb{C}(E)\tensor \Lambda^{(1,0)}(C))$.
Hitchin \cite{H1} proposes to add another equation to reduce the 
dimensionality. The \emph{Hitchin equations}
are a system of equations 
\begin{equation}
\label{eq:Hitchineq}
\begin{cases}
\overline{\partial}_A \phi = 0\\
F_A + [\phi,\phi^*] = 0 \,,
\end{cases}
\end{equation}
where 
$
d_A = \partial_A + \overline{\partial}_A
$
is the decomposition of the covariant derivative of 
the connection $A$ into its type $(1,0)$ and $(0,1)$ 
components that are determined by the complex structure of
$C$. The origin of (\ref{eq:Hitchineq}) is the dimensional 
reduction of the $4$-dimensional Yang-Mills theory. 
Hitchin observes that the 
 self-duality equation on $\mathbb{R}^4$
 restricted to $2$ dimensions by imposing 
independence in two variables automatically reduces to
(\ref{eq:Hitchineq}).

Since $A$ is a unitary connection in $E$, it defines a
holomorphic structure in $E$ through the covariant
Cauchy-Riemann operator $\overline{\partial}_A$. 
With respect to this complex structure, the first equation
 $\overline{\partial}_A \phi = 0$ implies that $\phi
 \in \Gamma(C,ad_\mathbb{C}(E)\tensor \Lambda^{(1,0)}(C))$
is  holomorphic. We recall that the holomorphic
part of $ad(E)$  is the holomorphic
endomorphism sheaf $\End(E)$ on $C$, and  the holomorphic
part of $\Lambda^{(1,0)}(C)$ is the sheaf of holomorphic
$1$-forms on $C$, or the \emph{canonical sheaf} $K_C$ on $C$.
Therefore, a solution of 
 $\overline{\partial}_A \phi = 0$ is a section
\begin{equation}
\label{eq:holomorphicphi}
\phi\in H^0(C,\End(E)\tensor K_C).
\end{equation}

We cannot define  the symplectic
quotient $\mathcal{A}(C,GL_n(\mathbb{C}))/\!\!/
\mathcal{G}(C,GL_n(\mathbb{C}))$ directly as we did before,
because $GL_n(\mathbb{C})$ is not compact and the 
analysis we need to deal with the infinite-dimensional
manifolds does not work.
The argument of Atiyah and Bott we have used in 
Section~\ref{sect:twisted} can be certainly applied to 
$\rho\in \mathcal{X}(\mathbb{C})$ of (\ref{eq:X}),
resulting in a projectively flat $gl_n(\mathbb{C})$
Yang-Mills connection $A_\mathbb{C}$ on $C$. 
It's $(0,1)$ part defines a holomorphic 
structure in the topological vector bundle $E$ as before,
but since the connection is not unitary, 
we are utilizing  only half of the information that 
$A_\mathbb{C}$ has. Hitchin's idea is that the other
half of the information goes to 
$\phi\in H^0(C,\End(E)\tensor K_C)$ through the 
factorization  (\ref{eq:Aphi}). Now the Serre duality
$$
H^0(C,\End(E)\tensor K_C)
= H^1(C,\End(E))^*
$$
and the Kodaira-Spencer deformation theory 
$$
H^1(C,\End(E)) = T_E \mathcal{U}_C(n,d)
$$
show that the pair $(E,\phi)$
is indeed 
an element of  $T^*\mathcal{U}_C(n,d)$,
which is what we expected in (\ref{eq:cotangentinX}). 
This pair consisting of  a holomorphic
vector bundle $E$ and a \emph{Higgs field}
$\phi$ of (\ref{eq:holomorphicphi}) is  known as
a \emph{Higgs pair} or a \emph{Higgs bundle}.

There is a slight inaccuracy here because we did not
impose any stability condition on $E$. The right notion
of stability
is that the slope inequality (\ref{eq:semistable}) holds
for every $\phi$-invariant proper vector subbundle $F$. 
Then the moduli space of unitary gauge equivalent
classes of irreducible solutions of the Hitchin equations 
(\ref{eq:Hitchineq}) is diffeomorphic to the moduli
space of stable Higgs pairs. Here we are assuming 
that the rank and the degree of $E$ are relatively prime. 
Obviously, if $E$ itself is stable, then the Higgs bundle
$(E,\phi)$ is stable for every $\phi$ in
$H^0(C,\End(E)\tensor K_C)$. Therefore, the complex
cotangent bundle  $T^*\mathcal{U}_C(n,d)$ is contained
in the moduli space $\mathcal{H}_C(n,d)$
of stable Higgs bundles as an open dense
subset. We also note that  the stability of a 
Higgs pair $(E,0)$ simply means that $E$ is stable.

Now we come back to the action 
of the  group
$\mathcal{G}(C,GL_n(\mathbb{C}))$ 
of complex gauge transformation on 
the space of complex valued connections
$\mathcal{A}(C,GL_n(\mathbb{C}))$. 
As we have noted earlier, we cannot directly define the symplectic 
quotient. After reducing the problem to 
considering Higgs pairs $(E,\phi)$, still we have
the ambiguity of the action of $H^0(C,\Aut(E))$ on 
the pairs since $E$ is not necessarily stable. 
But this situation is better than the symplectic quotient,
because of the fact that for every \emph{stable} Higgs pair $(E,\phi)$,
we have \cite{H1}
$$
H^0(C,\End(E,\phi)) = \mathbb{C}.
$$
Here an endomorphism of a Higgs bundle $(E,\phi)$ is defined to be a holomorphic
endomorphism $\psi$ of $E$ that commutes with $\phi$:
\begin{equation*}
\begin{CD}
E @>\psi>> E\\
@V{\phi}VV @VV{\phi}V\\
E\tensor K_C@>>{\psi\tensor 1}> E\tensor K_C 
\end{CD}
\end{equation*}

Although we know topological structures such
as the Poincar\'e polynomial of $T^*\mathcal{U}_C(n,d)$
from the work of \cite{AB} and \cite{HN},
their methods do not directly apply to the study  of the
character 
variety $\mathcal{X}(\mathbb{C})$. 
The work of Hausel and his collaborators \cite{HRVK}
reveals unexpectedly rich structures in the study of
the topology of these complex character varieties,
such as an unexpected relation to Macdonald polynomials.

\section{Hitchin integrable systems}
\label{sect:Hitchin}

From the point of view of 
$2$-dimensional Yang-Mills theory, we 
are led to identifying  the complex character variety
$\mathcal{X}(\mathbb{C})$ as the moduli
space $\mathcal{H}_C(n,d)$ of stable Higgs bundles.
In this section we show that there is an algebraically
completely integrable system on 
this Hitchin moduli space.

 The total space of the
complex cotangent bundle 
$T^*\mathcal{U}_C(n,d)$   is an open non-singular
complex submanifold of 
 $\mathcal{H}_C(n,d)$. 
 Since the cotangent bundle is easier to understand than
the Hitchin moduli, let us look at it first.    Note that
 $p^* \Lambda^1 (\mathcal{U}_C(n,d))
\subset \Lambda^1(T^*\mathcal{U}_C(n,d))$
has a tautological section 
$$
\eta\in H^0(T^*\mathcal{U}_C(n,d), 
p^* \Lambda^1 (\mathcal{U}_C(n,d))),
$$
where $
p:T^*\mathcal{U}_C(n,d)\rightarrow 
\mathcal{U}_C(n,d)
$ is the projection, and 
 $\Lambda^r (X)$ denotes in this section the sheaf of
holomorphic $r$-forms on a complex manifold $X$. 
The differential $\omega_I = d\eta$ of the 
tautological section defines the canonical holomorphic symplectic
form  on $T^*\mathcal{U}_C(n,d)$. The suffix $I$ indicates
the referrence to the complex structure of $\mathcal{U}_C(n,d)$.
The restriction of $\omega_I$ on $\mathcal{U}_C(n,d)$,
which is the $0$-section of the cotangent bundle, is identically
$0$. Therefore the $0$-section is a Lagrangian submanifold
of this holomorphic symplectic manifold. 

A surprising result of 
another influential paper \cite{H2} of Hitchin's is that
$\mathcal{H}_C(n,d)$ is the
total space of a  Lagrangian torus fibration. 
The starting point of his discovery is the following
intriguing equality as a consequence of the Riemann-Roch formula:
$$
\dim_C \mathcal{U}_C(n,d) = n^2(g-1) + 1
=1+  (g - 1) \sum_{i=1} ^ n (2 i-1) = \dim_C \bigoplus_{i=1} ^ n 
H^0(C,K_C ^{\tensor i}).
$$
Let us denote by
\begin{equation}
\label{eq:baseGL}
V_{GL} = V_{GL_n(\mathbb{C})} = 
\bigoplus_{i=1} ^ n 
H^0(C,K_C ^{\tensor i}).
\end{equation}
As a vector space $V_{GL}$ has the same 
dimension as 
$
H^0(C,\End(E)\tensor K_C) = T^* _E \,\mathcal{U}_C(n,d)
$.
The Higgs field $\phi\in H^0(C,\End(E)\tensor K_C)$
introduced by Hitchin earlier in \cite{H1} is a ``twisted'' endomorphism
$$
\phi: E \longrightarrow E\tensor K_C ,
$$
which induces a homomorphism of the $i$-th anti-symmetric 
tensor product spaces
$$
\wedge^i (\phi): \wedge^i (E) \longrightarrow 
\wedge^i (E\tensor K_C) = \wedge^i(E) \tensor K_C ^{\tensor i},
$$
or equivalently $\wedge^i (\phi)\in
H^0(C,\End(\wedge^i(E))\tensor K_C ^{\tensor i})$.
Taking its natural trace, we obtain
$$
\tr \wedge^i(\phi) \in H^0(C,K_C ^{\tensor i}).
$$
This is exactly the $i$-th characteristic coefficient of
the twisted endomorphism $\phi$:
\begin{equation}
\label{eq:characteristic}
\det (x - \phi) = x^n + \sum_{i=1} ^n (-1)^i \tr \wedge^i(\phi) \cdot
x^{n-i}.
\end{equation}
By assigning its coefficients, Hitchin \cite{H2} defines a
holomorphic map, now known as the \emph{Hitchin fibration}
or \emph{Hitchin map},
\begin{equation}
\label{eq:Hitchinmap}
H: \mathcal{H}_C (n,d)\owns (E,\phi)\longmapsto
\det(x-\phi)\in \bigoplus_{i=1} ^ n 
H^0(C,K_C ^{\tensor i}) = V_{GL}.
\end{equation}
The map $H$ to a vector space $V_{GL}$ is a collection
of $N= n^2(g-1)+1$ globally defined holomorphic
functions on $\mathcal{H}_C (n,d)$. The  $0$-fiber of the Hitchin fibration is the moduli space $\mathcal{U}_C(n,d)$.

What are other fibers of $H$? To answer this question,
the notion of \emph{spectral curves} is introduced in 
\cite{H2}. Generically other fibers  are
the Jacobians of these spectral curves. The total space of the canonical sheaf $K_C = \Lambda^1(C)$
on $C$ is the cotangent bundle $T^*C$. Let 
$$
\pi:T^*C \longrightarrow C
$$
be the projection, and 
$$
\tau\in H^0(T^*C, \pi^*K_C)\subset H^0(T^*C,\Lambda^1(T^*C))
$$
be the tautological section of $\pi^*K_C$ on $T^*C$. Here again
$\omega = d\tau$ is the holomorphic symplectic
form on $T^*C$. The tautological section $\tau$ satisfies
that $\sigma^*\tau = \sigma$
for every section $\sigma \in H^0(C,K_C)$ viewed as a holomorphic
map $\sigma:C\rightarrow T^*C$. The characteristic coefficients 
(\ref{eq:characteristic})
of $\phi$ 
 give a section
\begin{equation}
\label{eq:chareq}
s=\det(\tau-\phi) = \tau ^{\tensor n} + 
\sum_{i=1} ^n (-1)^{i}
\tr \wedge^i(\phi)\cdot \tau^{\tensor n-1}
\in H^0(T^*C,\pi^* K_C ^{\tensor n}).
\end{equation}
We define the spectral curve $C_s$ associated with a Higgs pair
$(E,\phi)$ as the divisor of $0$-points of the section $s=\det(\tau-\phi)$ of
the line bundle $\pi^*K_C ^{\tensor n}$:
\begin{equation}
\label{eq:spectralcurve}
C_s = (s)_0 \subset T^*C .
\end{equation}
 The 
spectral curve is the locus of $\tau$ that satisfies the 
characteristic equation $\det(\tau-\phi) = 0$. Thus 
every point of $C_s$  is
an eigenvalue, or spectrum,  of the twisted endomorphism $\phi$. This 
is the origin of the name of $C_s$. The projection $\pi$ defines
a ramified covering map 
$
\pi:C_s\rightarrow C
$ of degree $n$.

Another way to look at the spectral curve $C_s$ is to go 
through algebra. It has an advantage in identifying 
the fibers of the Hitchin fibration.
Since the section $s=\det(\tau-\phi)$ is
determined by the characteristic coefficients of $\phi$, by
abuse of notation we consider $s$ as an element of $V_{GL}$:
$$
s = (s_1,s_2,\dots, s_n) = (-\tr\, \phi,
\tr \wedge^2 (\phi), \dots, (-1)^n \tr \wedge^n(\phi))
\in \bigoplus_{i=1} ^n H^0(C,K_C ^{\tensor i}).
$$
It defines an $\mathcal{O}_C$-module 
$(s_1+s_2+\cdots +s_n)\tensor K_C ^{\tensor -n}$. Let 
$\mathcal{I}_s$ denote  the ideal generated by this module
in the symmetric tensor algebra $\Sym(K_C ^{-1})$. Since
$K_C ^{-1}$ is the sheaf of linear functions on $T^*C$, the scheme
associated to this tensor algebra is 
$\Spec \big(\Sym (K_C ^{-1}) \big)= T^*C$. The 
spectral curve as the divisor 
of $0$-points of $s$ is then defined by
\begin{equation}
\label{eq:spectralcurve2}
C_s = \Spec\left(\frac{\Sym(K_C ^{-1})}{\mathcal{I}_s}\right)
\subset \Spec \big(\Sym (K_C ^{-1}) \big)= T^*C .
\end{equation}
The set $U$
 consisting of points $s$ for which $C_s$ is irreducible 
 and non-singular
is an open dense subset
of $V_{GL}$ \cite{BNR}. The genus of $C_s$ can be found
as follows. Note that we have
$$
\pi_* \mathcal{O}_{C_s} = \Sym (K_C ^{-1})/\mathcal{I}_s
\isom \bigoplus_{i=0} ^{n-1} K_C^{\tensor -i}
$$ 
as an $\mathcal{O}_C$-module. From  the Riemann-Roch formula
we see that
$$
1-g(C_s) = \rchi(C_s,\mathcal{O}_{C_s})
= \rchi(C, \pi_* \mathcal{O}_{C_s}) = (1-g(C))\sum_{i=0} ^{n-1}
(2i+1)
= n^2(1- g(C)).
$$
Hence $g(C_s) = n^2(g-1) + 1$. As a consequence, 
we notice that the dimensions of the Jacobian variety
$\Jac(C_s)$ and the moduli space $\mathcal{U}_C(n,d)$ are
the same. The theory of spectral curves \cite{BNR, H2} 
makes this equality into  a precise geometric relation between these 
two spaces.

The Higgs field $\phi\in H^0(C,\End(E)\tensor K_C)$ gives
a homomorphism 
$$
\varphi: K_C ^{-1} \longrightarrow \End(E),
$$
which induces an algebra homomorphism, still denoted by the
same letter,
$$
\varphi : \Sym (K_C ^{-1}) \longrightarrow \End(E).
$$
Thus $\varphi$ defines a  $\Sym (K_C ^{-1})$-module
structure in $E$. 
Since $s\in V_{GL}$ is the characteristic coefficients of $\varphi$,
by the Cayley-Hamilton theorem,
 the homomorphism $\varphi$ factors through
$$
\Sym (K_C ^{-1}) \longrightarrow \Sym (K_C ^{-1})/\mathcal{I}_s
\longrightarrow
\End(E).
$$
Hence $E$ is actually a module over 
$\Sym (K_C ^{-1})/\mathcal{I}_s$ of rank $1$. The rank is 
$1$ because
the ranks of $E$ and
 $\Sym (K_C ^{-1})/\mathcal{I}_s$
 are the same as $\mathcal{O}_C$-modules. 
In this way a Higgs pair $(E,\phi)$ gives rise to a line bundle
$\mathcal{L}_E$ on the spectral curve $C_s$, if it is non-singluar.
Since $\mathcal{L}_E$ being an $\mathcal{O}_{C_s}$-module
is equivalent to $E$ being a 
$\Sym (K_C ^{-1})/\mathcal{I}_s$-module, we recover $E$
from $\mathcal{L}_E$ simply by
$E = \pi_* \mathcal{L}_E$, which has rank $n$ because $\pi$ is a covering
of degree $n$. From the equality 
$\rchi(C,E) = \rchi(C_s,\mathcal{L}_E)$ and Riemann-Roch, we find
that $\deg \mathcal{L}_E = d+ n(n-1) (g-1)$.
To summarize, the above construction defines an inclusion
map
$$
H^{-1}(s) \subset  \Pic^{d+n(n-1)(g-1)}(C_s) \isom \Jac(C_s),
$$
if $C_s$ is irreducible and non-singular.

Conversely,
suppose we have a line bundle $\mathcal{L}$ of degree $d+n(n-1)(g-1)$
on an irreducible non-singular spectral curve 
$C_s$. Then $\pi_*\mathcal{L}$ is a module over
$\pi_*\mathcal{O}_{C_s}
=\Sym (K_C ^{-1})/\mathcal{I}_s$, which
defines a homomorphism $\psi:K_C^{-1}\rightarrow
\End(\pi_*\mathcal{L})$. It is easy to 
see that the Higgs pair $(\pi_*\mathcal{L},\psi)$ is stable. Suppose
we had a $\psi$-invariant subbundle $F\subset \pi_*\mathcal{L}$
of rank $k<n$. Since $(F,\psi|_F)$ is a Higgs pair,
it gives rise to a spectral curve $C_{s'}$. From the construction,
we have an injective morphism $C_{s'}\rightarrow C_s$. But since
$C_s$ is irreducible, it contains no smaller component.
Therefore, $\pi_*\mathcal{L}$ has no $\psi$-invariant proper subbundle. 
Thus we have established that
\begin{equation}
\label{eq:Hitchinfiber}
H^{-1}(s) \isom \Jac(C_s), \qquad s\in U\subset V_{GL}.
\end{equation}
We note that the vector bundle $\pi_*\mathcal{L}$ is not necessarily 
stable. 
It is proved in \cite{BNR} that the locus of $\mathcal{L}$ in
$\Pic^{d+n(n-1)(g-1)}(C_s)$ that gives unstable
$\pi_*\mathcal{L}$ has codimension two or more.

Recall that  the tautological section
$\eta\in H^0(T^*\mathcal{U}_C(n,d),p^*\Lambda^1(
\mathcal{U}_C(n,d)))$ is a  holomorphic
$1$-form on $T^*\mathcal{U}_C(n,d)
\subset \mathcal{H}_C(n,d)$. Its restriction to the fiber $H^{-1}(s)$
of  $s\in U$ for which $C_s$ is non-singular
extends to a holomorphic $1$-form on the whole 
fiber $H^{-1}(s)\isom \Jac(C_s)$ since
$\eta$ is undefined only on a codimension $2$ subset.
Consequently $\eta$ extends as a holomorphic $1$-form
on $H^{-1}(U)$.
Thus $\eta$ is well defined on 
$T^*\mathcal{U}_C(n,d)\cup H^{-1}(U)$. The complement of
this open subset in $\mathcal{H}_C(n,d)$ consists of such Higgs pairs
$(E,\phi)$ that $E$ is unstable \emph{and} $C_s$ is 
singular. Since the stability of $E$ and the non-singular condition 
for $C_s$ are
both open conditions, this complement has codimension at least 
two. Consequently, both the tautological section $\eta$ and the 
holomorphic symplectic form $\omega_I = d\eta$ extend holomorphically
to the whole Higgs moduli space $\mathcal{H}_C(n,d)$. 

We note that there are no holomorphic $1$-forms other than  constants
on a Jacobian variety.  It implies  that
$$
\omega_I|_{H^{-1}(s)} = d(\eta|_{H^{-1}(s)}) = 0
$$
for $s\in U$. 
The \emph{Poisson structure}
on $H^0(\mathcal{H}_C(n,d),\mathcal{O}_{\mathcal{H}_C(n,d)})$ is defined by
$$
\{f,g\} = \omega_I(X_f,X_g),\qquad f,g\in 
H^0(\mathcal{H}_C(n,d),\mathcal{O}_{\mathcal{H}_C(n,d)}),
$$
where $X_f$ denotes the Hamiltonian vector field defined by
 the relation $df= \omega_I(X_f,\cdot)$. 
Since $\omega_I$ vanishes on a generic  fiber of $H$, the holomorphic functions on 
$\mathcal{H}_C(n,d)$ coming from coordinates
of the Hitchin fibration
 are \emph{Poisson commutative} with respect to the
holomorphic symplectic structure $\omega_I$. 
An \emph{algebraically completely integrable Hamiltonian system} on a holomorphic symplectic manifold $(M, \omega)$
of dimension $2m$  is an 
open
 holomorphic
map $H:M\rightarrow \mathbb{C}^m$ such that the coordinate functions
are Poisson commutative and a generic fiber is an Abelian variety
\cite{V}. Thus $(\mathcal{H}_C(n,d),\omega_I, H)$ is an algebraically
completely integrable Hamiltonian system,
called the \emph{Hitchin integrable system}. 

\begin{thm}
The Hitchin fibration
$$
H: \mathcal{H}_C(n,d)\longrightarrow V_{GL}
$$
is a Lagrangian Jacobian fibration
defined on an algebraically completely integrable
system $(\mathcal{H}_C(n,d),\omega_I, H)$.  A
generic fiber $H^{-1}(s)$
is a Lagrangian with respect to the holomorphic symplectic 
structure $\omega_I$ and is isomorphic to the Jacobian variety
of a  spectral curve $C_s$.
\end{thm}

\section{Symplectic quotient of the Hitchin
system and mirror symmetry}
\label{sect:quotient}

Is the Hitchin fibration (\ref{eq:Hitchinmap}) an
effective family of deformations of Jacobians?
This is the question we address
in \cite{HM}. The investigation of this question leads to the
relation between the Hitchin systems and mirror symmetry
discovered by Hausel and Thaddeus \cite{HT}.

The Jacobian variety $\Jac(C)=\Pic^0(C)$ acts on $\mathcal{H}_C(n,d)$
by 
$
(E,\phi)\longmapsto (E\tensor L,\phi)$, where $L\in \Jac(C)$
is a line bundle on $C$ of degree $0$. The Higgs field is preserved
because 
$$
E^*\tensor E \longmapsto (E\tensor L)^*\tensor (E\tensor L)=E^*
\tensor E
$$
 is unchanged. Thus this action does not contribute to 
 deformations of the spectral curves. It is natural to
 symplectically quotient it out.  On the open subset 
$T^*\mathcal{U}_C(n,d)$,  the $\Jac(C)$ action is symplectomorphic
because it is
induced by the action on the base space $\mathcal{U}_C(n,d)$.
On the other open subset $H^{-1}(U)$ the action is also symplectomorphic 
because it preserves each fiber which is a Lagrangian. Thus 
 the action of $\Jac(C)$ on 
$\mathcal{H}_C(n,d)$ is globally symplectomorphic. We claim
that the first component of the Hitchin map
$$
H_1:\mathcal{H}_C(n,d)\owns (E,\phi)\longmapsto
\tr(\phi) \in H^0(C,K_C)
$$
is the moment map of this Jacobian action. Note that 
$H^1(C,\mathcal{O}_C)$ is the Lie algebra of the Abelian
group $\Jac(C)$, hence $H^0(C,K_C)$ is the dual Lie algebra. 
The claim is obvious because $\omega_I$ vanishes on each
fiber of the Hitchin fibration on which the $\Jac(C)$ action is
restricted, and because $dH_1$ is the $0$-map on any
infinitesimal deformation of $E$. Therefore, we can define
the symplectic quotient
\begin{equation}
\label{eq:ph}
\mathcal{PH}_C(n,d)
\overset{\text{def}}{=}\mathcal{H}_C(n,d)/\!\!/\Jac(C) = H_1 ^{-1}(0)/\Jac(C).
\end{equation}
It's dimension is $2(n^2-1)(g-1)$. The letter P stands for
``projective.''

The moment map $H_1$ being the trace of $\phi$, it  is natural  to 
define
\begin{equation}
\label{eq:Vsl}
V_{SL} = V_{SL_n(\mathbb{C})} = 
\bigoplus_{i=2} ^n H^0(C,K_C ^{\tensor i}) \subset
V_{GL}.
\end{equation}
This is a vector space of dimension $(n^2-1)(g-1)$. Since
the $\Jac(C)$-action on $\mathcal{H}_C(n,d)$ preserves 
fibers of the Hitchin fibration, the map $H$ induces a natural
map
\begin{equation}
\label{eq:Hpgl}
H_{PGL}:\mathcal{PH}_C(n,d)\longrightarrow V_{SL}.
\end{equation}
It's $0$-fiber  is $H_{PGL} ^{-1}(0) = \mathcal{U}_C(n,d)/\Jac(C)$.
To study the symplectic quotient (\ref{eq:ph}),
 let us first analyze this $0$-fiber.
Following \cite{MFK} we denote by $\mathcal{SU}_C(n,d)$ the
moduli space of stable vector bundles with a fixed determinant line
bundle. This is a fiber of the determinant map
\begin{equation}
\label{eq:det}
\mathcal{U}_C(n,d)\owns E\longmapsto \det E\in \Pic^d(C),
\end{equation}
and is independent of the choice of the value of the determinant. 
This fibration  is a non-trivial fiber bundle. The
equivariant  $\Jac(C)$-action on (\ref{eq:det}) is given by
\begin{equation}
\label{eq:equivariantaction}
\begin{CD}
\mathcal{U}_C(n,d) @>{\tensor L}>> \mathcal{U}_C(n,d)\\
@V{\text{det}}VV @VV{\text{det}}V\\
\Pic^d(C) @>>{\tensor L^{\tensor n}}> \Pic^d(C)
\end{CD}
\qquad L\in \Jac(C).
\end{equation}
The isotropy subgroup of the $\Jac(C)$-action on 
$\Pic^d(C)$ is 
the group of $n$-torsion points
$$
J_n(C) \overset{\text{def}}{=}\{L\in \Jac(C)\,|\,
L^{\tensor n}=\mathcal{O}_C\} \isom
H^1(C,\mathbb{Z}/n\mathbb{Z}) .
$$
Choose a reference line bundle $L_0\in \Pic^d(C)$ and 
consider a degree $n$ covering 
$$
\nu : \Pic^d(C)\owns L\tensor L_0\longmapsto 
L^{\tensor n}\tensor L_0\in\Pic^d(C) , \qquad L\in \Jac(C) .
$$
Then the pull-back bundle $\nu^*\mathcal{U}_C(n,d)$
on $\Pic^d(C)$ becomes trivial:
$$
\nu^*\mathcal{U}_C(n,d) = \Pic^d(C)\times \mathcal{SU}_C(n,d) .
$$
The quotient of this product by the diagonal action of 
$J_n(C)$ is the original moduli space:
\begin{equation}
\label{eq:Jnquotient1}
\big( \Pic^d(C)\times \mathcal{SU}_C(n,d)\big)\big/
J_n(C) \isom \mathcal{U}_C(n,d).
\end{equation}
It is now clear that 
$$
\mathcal{U}_C(n,d)/\Jac(C) \isom \mathcal{SU}_C(n,d)/
J_n(C) .
$$

The other fibers of (\ref{eq:Hpgl}) are best described in terms of 
\emph{Prym varieties}. Let $s\in V_{SL}\cap U$ be a point such that
$C_s$ is irreducible and non-singular. 
The covering map 
$\pi:C_s \rightarrow C$ induces an injective  homomorphism
$
\pi^*: \Jac(C)\owns L \longmapsto \pi^*L\in  \Jac(C_s)
$.
This is injective because if $\pi^*L\isom \mathcal{O}_{C_s}$,
then by the projection formula we have
$$
\pi_*(\pi^*L) \isom \pi_*\mathcal{O}_{C_s}\tensor L\isom
\bigoplus_{i=0} ^{n-1} L\tensor K_C ^{\tensor -i},
$$
which has a nowhere vanishing section.
Hence $L\isom \mathcal{O}_C$. Take a point $(E,\phi)\in
H^{-1}(s)$ and let $\mathcal{L}_E$ be the corresponding line 
bundle on $C_s$. Since 
$\pi_*(\mathcal{L}_E\tensor \pi^*L) \isom E\tensor L$, 
the action of $\Jac(C)$ on $H^{-1}(s)\isom 
\Jac(C_s)$ is the canonical subgroup action. Thus we conclude that 
the fiber $H_{PGL} ^{-1}(s)$ is isomorphic to the 
\emph{dual Prym variety} of the covering $C_s\rightarrow C$
\begin{equation}
\label{eq:dualprym}
{\Prym}^*(C_s/ C) \overset{\text{def}}{=} \Jac(C_s)/\Jac(C).
\end{equation}
The \emph{Prym variety} ${\Prym}(C_s/ C)$ 
of the covering is defined to be 
the kernel of the norm map
\begin{equation}
\label{eq:norm}
\Nm: \Jac(C_s)\owns \mathcal{L}\longmapsto
\det(\pi_*\mathcal{L})\tensor (\det\pi_*\mathcal{O}_{C_s})^*
\in \Jac(C).
\end{equation}
Both Prym and dual Prym varieties are Abelian varieties of
dimension $g(C_s)-g(C)$. 
Similarly to the equivariant action (\ref{eq:equivariantaction}),
we have 
\begin{equation}
\label{eq:equivariantjac}
\begin{CD}
\Jac(C_s) @>{\tensor L}>> \Jac(C_s)\\
@V{\Nm}VV @VV{\Nm}V\\
\Jac(C) @>>{\tensor L^{\tensor n}}> \Jac(C)
\end{CD}
\qquad L\in \Jac(C).
\end{equation}
By the same argument as we used in (\ref{eq:Jnquotient1}), we obtain
\begin{equation}
\label{eq:Jnquotient2}
\big(\Prym(C_s/C)\times \Jac(C)\big)\big/J_n(C) \isom \Jac(C_s).
\end{equation}
From (\ref{eq:dualprym}) and (\ref{eq:Jnquotient2}), it follows
that
$
{\Prym}^*(C_s/ C) =\Prym(C_s/ C) /J_n(C)
$.
We have thus established

\begin{thm}
The fibration
$H_{PGL}:\mathcal{PH}_C(n,d)\rightarrow V_{SL}$ is
a generically Lagrangian dual Prym fibration.
\end{thm}

How can we construct a Lagrangian \emph{Prym} fibration?
The dual Prym variety naturally appears in the above discussion
when we quotient out the Jacobian action on the moduli space
of vector bundles.
Another way to limit the Jacobian action is to restrict the
structure group of the vector bundles 
from $GL_n(\mathbb{C})$ to $SL_n(\mathbb{C})$. So let us
consider a character variety
$$
\Hom
 (\hat{\pi}_1(C)_\mathbb{R}, SL_n(\mathbb{C}))/\!\!/
SL_n(\mathbb{C}).
$$
Although the central generator $c\in \hat{\pi}_1(C)$ can take the 
same value as in (\ref{eq:gammavalue}), to have a representation of
$\hat{\pi}_1(C)_\mathbb{R}$, $c$ has to be mapped to the identity. 
Thus we go back to the untwisted character variety
$\Hom
 ({\pi}_1(C), SL_n(\mathbb{C}))/\!\!/
SL_n(\mathbb{C})$. 
The argument of Section~\ref{sect:gln} leads us to 
the moduli space of stable Higgs pairs $(E,\phi)$,
where this time $\det(E)=\mathcal{O}_C$ and the 
Higgs field $\phi:E\rightarrow E\tensor K_C$
is traceless 
since $\End(E)$ is an $sl_n(\mathbb{C})$-bundle.
 Let us denote this moduli space by
$\mathcal{SH}_C(n,0)$. 
Here the letter S stands for 
``special.'' 
The natural counterpart of
 the Hitchin fibration  on $\mathcal{SH}_C(n,0)$
is the map
\begin{equation}
\label{eq:Hsl}
H_{SL}: \mathcal{SH}_C(n,0)\owns (E,\phi)
\longmapsto \det(x-\phi)\in V_{SL}.
\end{equation}
It's $0$-fiber is $H_{SL} ^{-1}(0) = \mathcal{SU}_C(n,0)$. 
For a generic $s\in V_{SL}$ for which $C_s$ is irreducible and
non-singular, obviously we have $H_{SL} ^{-1}(s) \isom
\Prym(C_s/C)$. 

\begin{thm}[ \cite{HT, DP}] The two
Lagrangian Abelian fibrations 
\begin{equation}
\label{eq:mirrors}
\begin{CD}
\mathcal{SH}_C(n,0) @.  \mathcal{PH}_C(n,d)\\
@V{H_{SL}}VV @VV{H_{PGL}}V\\
V_{SL} @= V_{SL} 
\end{CD}
\end{equation}
are \textbf{mirror dual} in the sense of Strominger-Yau-Zaslow
\cite{SYZ}. 
\end{thm}
The mirror duality here means that the 
bounded derived category $D^b(Coh(\mathcal{SH}_C(n,0)))$
of coherent analytic sheaves on
$\mathcal{SH}_C(n,0)$
is equivalent to the Fukaya category 
$Fuk(\mathcal{PH}_C(n,d))$ consisting of Lagrangian 
subvarieties 
of $\mathcal{PH}_C(n,d)$ and flat $U_1$-bundles on them
\cite{Fukaya}. We can view it as a family of deformations
 of Furier-Mukai
duality \cite{Mukai, P} between 
$\Prym(C_s/C)$ and $\Prym^*(C_s/C)$
 parametrised 
on the same base space $V_{SL}$.

As noted at the end of Section~\ref{sect:bundle}, $\Jac(C)$ of
an algebraic curve $C$ is the moduli space of flat $U_1$ 
connections modulo gauge transformation. This correspondence
does not require that $C$ is a curve, because  the
flatness condition automatically implies the
integrability of the $(0,1)$-part of the connection. 
Since the Abel-Jacobi
map $C\rightarrow \Jac(C)$ induces a homology isomorphism
$$
H_1(C,\mathbb{Z})\overset{\sim}{\longrightarrow}
H_1(\Jac(C),\mathbb{Z}),
$$
we have an isomorphism
$$
\Pic^0(\Jac(C))\overset{\sim}{\longrightarrow}\Jac(C),
$$
because any representation of the fundamental group in
$U_1$ factors through the Abelian group homomorphism
from the homology group.
Here $\Pic^0$ indicates the moduli of holomorphic
line bundles that are topologically trivial.
Thus $\Jac(C)$ is \emph{self-dual}. Now consider a flat $U_1$ connection $A$ on 
$\Prym^*(C_s/C)$. It is  a holomorphic 
line bundle on $\Jac(C_s)$ that is invariant under the $\Jac(C)$-action. 
The restriction of $A$ to $C\subset \Jac(C)\subset\Jac(C_s)$ 
then defines a holomorphic line bundle on $C$, which is trivial
by the assumption. We notice that this correspondence
$\Jac(C_s)\rightarrow \Jac(C)$ is exactly the norm map of
(\ref{eq:norm}). In other words, we obtain the duality
\begin{equation}
\label{eq:prymduality}
\Pic^0(\Prym^*(C_s/C)) \isom \Prym(C_s/C).
\end{equation}
A skyscraper sheaf on $\mathcal{SH}_C(n,0)$
supported at a point $(E,\phi)$
 determines a spectral curve
$C_s$ and a point on the Prym variety 
$\Prym(C_s/C)$, where
$s=H_{SL}(E,\phi)$. It then identifies a fiber
$H_{PGL} ^{-1}(s)\isom \Prym^*(C_s/C)$,
which is a Lagrangian subvariety of $\mathcal{PH}_C(n,d)$,
and a flat $U_1$-connection on it
because of (\ref{eq:prymduality}). This is the idea of geometric
realization of mirror symmetry due to Strominger, Yau and
Zaslow \cite{SYZ}. 

Although complex structures are different, we can identify
\begin{equation}
\begin{cases}
\mathcal{SH}_C(n,0) \isom \Hom(\pi_1(C),SL_n(\mathbb{C}))
/\!\!/SL_n(\mathbb{C}) \\
\mathcal{PH}_C(n,0) \isom \Hom(\pi_1(C),PGL_n(\mathbb{C}))
/\!\!/PGL_n(\mathbb{C}).
\end{cases}
\end{equation}
Then  the mirror symmetry  (\ref{eq:mirrors}) gives a manifestation of
\emph{geometric Langlands correspondence} \cite{DP, HT, KW},
which is a family of Fourier-Mukai duality  transformations over
 the same base
space \cite{Frenkel}. Thus \emph{the Hitchin integrable systems 
on character varieties relate
the SYZ mirror symmetry and the geometric Langlands correspondence}.

We have noted earlier that $\mathcal{H}_C(n,d)$ has two
different complex structures $I$ and $J$.
The complex structure $I$ comes from the moduli space of stable
Higgs bundles, and $J$ from a connected component 
$\mathcal{X}(\mathbb{C})$ of 
the twisted character variety
$\Hom(\hat{\pi}_1(C)_\mathbb{R},GL_n(\mathbb{C}))
/\!\!/GL_n(\mathbb{C})$.
The complex manifold $\mathcal{U}_C(n,d)$, assuming
$G.C.D.(n,d) = 1$, is projective algebraic, hence has a unique
K\"ahler metric. The K\"ahler form in a real coordinate is a real 
symplectic form, which extends to a holomorphic symplectic 
form $\omega_J$ on the complexification 
$\mathcal{X}(\mathbb{C})$ of $\mathcal{U}_C(n,d)$.
Thus $\omega_J ^N$ defines
a holomorphic top form on $\mathcal{X}(\mathbb{C})$,
where $N={\dim_\mathbb{C} \mathcal{U}_C(n,d)}$.
We can then think of $(\mathcal{X}(\mathbb{C}), J, \omega_J ^N,
\omega_I)$ as a $2N$-dimensional Calabi-Yau manifold. 
The Hitchin fibration is an example of a \emph{special Lagrangian
fibration}, meaning that the restriction of $\omega_J ^N$ on
each fiber $H^{-1}(s)$ gives a Riemannian volume form on 
$\Jac(C_s)$. Since 
$$
p:H^{-1}(s)\isom \Jac(C_s)\longrightarrow \mathcal{U}_C(n,d)
$$
is a finite covering of degree $2^{3(g-1)}\cdot 3^{5(g-1)}\cdots
n^{(2n-1)(g-1)}$ \cite{BNR}, a generic fiber $H^{-1}(s)$ has 
the same Riemannian volume 
that is equal to $2^{3(g-1)}\cdot 3^{5(g-1)}\cdots
n^{(2n-1)(g-1)}$-times the K\"ahler volume  of
$\mathcal{U}_C(n,d)$. 
Actually, the space $\mathcal{H}_C(n,d)=\mathcal{X}(\mathbb{C})$
is a \emph{hyper K\"ahler manifold} with complex structures 
$I$, $J$, and  $K=IJ$.

Kapustin and Witten \cite{KW} noticed that the mirror symmetry
(\ref{eq:mirrors}) is a consequence of the dimensional reduction 
of $4$-dimensional super Yang-Mills theory. In their formulation, the
Langlands duality corresponds to the physical electro-magnetic 
duality, and the Fourier-Mukai transform on each fiber of the
Hitchin fibrations is the $T$-duality.

Finally, let us comment on 
the relation between the Hitchin systems, Prym 
varieties, and Sato Grassmannians established in
 \cite{HM, LM1, LM2}.  A theorem of \cite{LM2}
basically  states that  to any morphism
$\pi:C_s\rightarrow C$ between algebraic curves, a solution
of a KP-type integrable system (the $n$-component KP
equations and more general Heisenberg KP equations)
is constructed with the following
two properties: a) the 
orbit of the dynamical system 
on the Grassmannian is the Prym variety $\Prym(C_s/C)$; and b)
the evolution equations are linearlized  on the Prym variety. 
To make a connection between the Hitchin integrable
systems and the theory of \cite{LM2}, we need to quotient out
the trivial 
deformations of spectral curves $\{C_s\}_{s\in V_{GL}}$
given by a scalar action
\begin{equation}
\label{eq:weighted}
V_{GL} = \bigoplus_{i=1} ^n H^0(C,K_C^{\tensor i})\owns
(s_1,s_2,\dots,s_n)\longmapsto
(\lam s_1,\lam^2 s_2,\dots, \lam^n s_n)\in V_{GL} 
\end{equation}
for $\lam\in \mathbb{C}^*$.
This action corresponds to  the
scalar multiplication of  
a Higgs field
$\phi\mapsto \lam\cdot\phi$, which  is \emph{not}
a symplectomorphism on the Hitchin moduli space 
because it changes the symplectic form
to $\lam\cdot\omega_I$. Let us define the projective
Hitchin moduli space
$$
\mathbb{P}(\mathcal{H}_C(n,d))
= \big(\mathcal{H}_C(n,d)\setminus H^{-1}(0)\big)\big/\mathbb{C}^* .
$$
This is no longer a holomorphic symplectic manifold, yet 
the Hitchin fibration naturally descends to a generically 
Jacobian fibration
$$
H_{GL} ^\mathbb{P}: \mathbb{P}(\mathcal{H}_C(n,d))\longrightarrow
\mathbb{P}_w(V_{GL})
$$
over the weighted projective space of $V_{GL}$ defined 
by (\ref{eq:weighted}). Now we have
\begin{thm}[\cite{HM}]
There is a rational map $\iota$ from $\mathbb{P}_w(V_{GL})$
into the Grassmannian of \cite{LM2} such that 
\begin{enumerate}
\item $\iota$ is generically an embedding;
\item the orbit of the $n$-component  KP equations 
starting at $\iota(\mathbb{P}_w(V_{GL}))$ in the 
Grassmannian 
is birational to $\mathbb{P}(\mathcal{H}_C(n,d))$; and
\item the dynamical system on $\mathbb{P}(\mathcal{H}_C(n,d))$
defined by the Hitchin integrable system is the pull-back of 
the $n$-component KP equations via $\iota$.
\end{enumerate}
\end{thm}

A similar theorem holds for the Prym fibration
(\ref{eq:Hsl}), where we use the traceless $n$-component 
KP equations to produce Prym varieties as orbits.

There is a common belief coming out of the recent developments 
on character varieties. It is
that 
to fully appreciate the categorical equivalences of the dualities
such as mirror symmetry and geometric Langlands correspondence,
the moduli theory based on stable objects is not the right language.
We are naturally led to considering \emph{moduli stack} of vector bundles
and other categorical objects. Infinite-dimensional geometry
of connections \cite{AB, H1} played an essential role in 
understanding the geometry of moduli spaces of stable vector bundles. 
Infinite-dimensional Sato Grassmannians are more suitable
geometric objects for algebraic stacks. Although our current
understanding of the relation between the Hitchin systems and
Sato Grassmannians is limited, 
more should be coming as our understanding of the duality 
deepens 
from this point of view.

\begin{ack}
I would like to thank the organizers of the RIMS-OCAMI
joint conference for their invitation and \emph{exceptional}
hospitality during my stay in Kyoto and Nara.
\end{ack}

\bigskip


\providecommand{\bysame}{\leavevmode\hbox to3em{\hrulefill}\thinspace}

\bibliographystyle{amsplain}

\end{document}